\documentclass[12pt]{article}
\usepackage{amssymb,amsmath,amsthm,epsfig,amscd,eucal}

\newtheorem{thm}{Theorem}[section]
\newtheorem{defn}[thm]{Definition}

\newcommand{\curlyG}{\mathcal{G}}

\title{Enumerating limit groups: A Corrigendum}
\author{Daniel Groves and Henry Wilton}
\date{\today}

\begin{document}
\maketitle

\begin{abstract}
We discuss two possible interpretations of Definition 1.1 from \cite{GW}.
\end{abstract}

In \cite[Definition 1.1]{GW}, we made the following definition:

\setcounter{section}{+1}

\begin{defn}\label{1.1}
A coherent group $G$ is \emph{effectively coherent} if there exists an algorithm that, given a finite subset $S$ as input, outputs a presentation for the subgroup generated by $S$.

A class $\curlyG$ of coherent groups is \emph{uniformly effectively coherent} if there exists an algorithm that, given as input a presentation of a group $G\in\curlyG$ and a finite set $S$ of elements of $G$, outputs a presentation for the subgroup of $G$ generated by $S$.
\end{defn}

We intended the phrase `presentation for the subgroup' to mean that the presentation has generating set $S$, so that one knows how the abstract group sits inside $G$ as a subgroup.  However, as pointed out to us by Maurice Chiodo, one could interpret our definition to mean that the algorithm merely outputs a presentation for the abstract group $\langle S \rangle$, without exhibiting an isomorphism between the group presented and $\langle S \rangle$.

In order to avoid further confusion, we propose the following definitions to distinguish the two related notions.

\begin{defn}\label{1.2}
A coherent group $G$ is \emph{effectively coherent} if there exists an algorithm that, given a finite subset $S$ that generates a subgroup $H$ as input, outputs a presentation $\langle S \mid R \rangle$ for $H$.
\end{defn}

We stress that the output presentation is required to be on the given input generating set $S$.  This is equivalent to requiring that we are given an isomorphism between $\langle S\rangle$ and the output presentation.

\begin{defn}\label{1.3}
A coherent group $G$ is \emph{weakly effectively coherent} if there exists an algorithm that, given a finite subset $S$ that generates a subgroup $H$ as input, outputs a presentation for $H$ as an abstract group.
\end{defn}

There are corresponding notions of \emph{uniform effective coherence} and \emph{uniform weak effective coherence} for classes of groups. 

We emphasise that, as long as Definition \ref{1.1} is interpreted as Definition \ref{1.2}, the results and proofs of \cite{GW} are correct as stated.  However, if one interprets it as Definition \ref{1.3} then some problems can arise---see \cite{Chiodo}, especially Theorem E.

Note that a locally Hopfian group (for example, a residually finite group) is effectively coherent if and only if it is weakly effectively coherent  (cf.\ \cite[Theorem F]{Chiodo}).
 
\medskip We thank Chiodo for pointing out the ambiguity in our definition, and apologise for the confusion.

\bibliographystyle{plain}

\noindent
\textsc{Daniel Groves, MSCS UIC 322 SEO, \textsc{M/C 249}, 851 S. Morgan St., Chicago, IL, 60607-7045, USA}\\
\emph{E-mail:} \texttt{groves@math.uic.edu}\\\\
\textsc{Henry Wilton, Department of Mathematics, University College London, Gower Street, London, WC1E 6BT, UK}\\
\emph{E-mail:} \texttt{hwilton@math.ucl.ac.uk}\\

\end{document}